\documentclass[11pt]{amsart}
\usepackage{graphics}
\usepackage{amsmath}
\usepackage{amsfonts}
\usepackage{latexsym}
\usepackage{enumerate}
\usepackage{array}
\usepackage{multirow}
\newcommand{\C}{{\bf C}}
\newcommand{\R}{{\bf R}}
\newcommand{\Q}{{\bf Q}}
\newcommand{\Z}{{\bf Z}}

\newcommand{\Aut}{{\rm Aut}}

\newcommand{\rtimes}{\times \kern -2.2pt
\vrule height 5.6pt depth 0pt width 0.37pt\;}

\begin{document}
\title{Complex and K\"ahler structures
on Compact Solvmanifolds
}
\author{Keizo Hasegawa\\
}
\begin{abstract}
We discuss our recent results on the existence and classification problem of
complex and K\"ahler structures on compact solvmanifolds. In particular,
we determine in this paper all the complex surfaces which are diffeomorphic
to compact solvmanifolds (and compact homogeneous manifolds in general).
\end{abstract}

\maketitle

\section{ Introduction}

The purpose of this paper is to discuss our recent results on the existence and
classification problem of complex and K\"ahler structures on compact solvmanifolds.
In particular, we determine in this paper all the complex surfaces which are diffeomorphic
to four-dimensional compact solvmanifolds (see Theorem 1); combined with
many known results, this makes us determine the complete list of
complex surfaces which are diffeomorphic to four-dimensional compact homogeneous
manifolds (see Section 8).
\smallskip

A solvmanifold $M$ is a compact homogeneous space of solvable Lie group,
that is, $M$ is a differentiable manifold on which a connected solvable Lie group $G$
acts transitively (and almost effectively). $M$ can be written as $D \backslash G$,
where $G$ is a simply connected solvable Lie group and $D$ is a closed subgroup
of $G$ (which includes no non-trivial connected normal subgroup of $G$). By complex
structures (or K\"ahler structures) on solvmanifolds we mean integrable almost complex
structures (with compatible K\"ahler form) on $M$ which are not necessarily invariant
by the canonical (right) action of $G$. A complex structure $J$ on $M$ is called
{\em left-invariant} if it is induced from a left-invariant complex structure on $G$. 
It should be noted that complex structures on solvmanifolds may or may not be
left-invariant;  however we shall see in this paper that four-dimensional solvmanifolds
admit only left-invariant complex structures. We do not know any
solvmanifolds (including tori) of higher dimension which admit
non-left-invariant complex structures.
\smallskip

The classification of compact homogeneous K\"ahler manifolds is well known
\cite{B, M}; in particular we know that the only compact homogeneous
K\"ahler solvmanifolds are complex tori.
On the other hand we observed in the paper \cite{H3} that a four-dimensional solvmanifold
admits a K\"ahler structure if and only if it is a complex torus or a hyperelliptic surface.
The following theorem on four-dimensional solvmanifolds may be considered as a
generalization of this result:
\medskip

\noindent {\bfseries Theorem 1.}  {\em A complex surface is diffeomorphic to a four-dimensional
solvmanifold if and only if it is one of the following surfaces:
Complex torus, Hyperelliptic surface, Inoue Surface of type $S^0$,
Primary Kodaira surface, Secondary Kodaira surface, Inoue Surface of type $S^{\pm}$.
And every complex structure on each of these complex surfaces (considered as solvmanifolds)
is left-invariant.}
\medskip

A hyperelliptic surface can be characterized as a finite quotient of a complex torus
which is simultaneously a complex torus bundle over a complex torus. As a natural
generalization of hyperelliptic surfaces to the higher dimension, we can define a class
of K\"ahlerian solvmanifolds (see Example 4). And as stated in the paper \cite{H3}, we made
a conjecture that any K\"ahlerian solvmanifold must belong to this class of solvmanifolds.
Recently we have given a complete proof for this conjecture, applying a result of
Arapura and Nori  on polycyclic (solvable) K\"ahler  groups \cite{AN}
(see also \cite{ABCKT}), together with our previous results on K\"ahlerian nilmanifolds
\cite{BG1, H1} and K\"ahlerian solvmanifolds \cite{H3} (see Section 7): 
\medskip

\noindent {\bfseries Theorem 2 {\rm (\cite{H4})}.} {\em A compact solvmanifold admits a
K\"ahler structure if and only if it is a finite quotient of a complex torus which
has a structure of a complex torus bundle over a complex torus. In particular,
a compact solvmanifold of completely solvable type has a K\"ahler structure
if and only if it is a complex torus.}
\medskip

The last part of the theorem was first conjectured by Benson and Gordon \cite{BG2},
which can be shown simultaneously in the proof of the theorem. For the definition of
{\em completely solvable type}, we refer to Section 7.
\medskip

In this paper we provide many examples of solvmanifolds, including two class of
four-dimensional solvmanifolds which admit no complex structures (see Section 4),
and a class of six-dimensional pseudo-K\"ahlerian solvmanifolds (of completely
solvable type) which admit no complex structures but satisfy most of the known
topological properties of compact K\"ahler manifolds (see Section 7).

\section{Fundamental results on solvmanifolds}

We recall some fundamental results on solvmanifolds, most of which are found
in \cite{AU2}. Let $M$ be a solvmanifold of dimension $n$. We have the following
basic results:

\begin{list}{}{\topsep=5pt \leftmargin=10pt \itemindent=10pt \itemsep=3pt}
\item[(1)] $M$ is a fiber bundle over a torus with fiber a nilmanifold (which is called
the {\em Mostow fibration} of $M$) \cite{MS}. In particular we can represent the
fundamental group $\Gamma$ of $M$ as an extension of a torsion-free
nilpotent group $\Lambda$ of rank $n-k$ by a free abelian group of rank $k$,
where $1 \le k \le n$ and $k = n$
if and only if $\Gamma$ is abelian:
$$0 \rightarrow \Lambda \rightarrow \Gamma \rightarrow \Z^k \rightarrow 0.$$

\item[(2)] Conversely any such abstract group $\Gamma$ (which is a polycyclic group)
can be the fundamental group of some solvmanifold \cite{WG}. We call such a
group $\Gamma$ a {\em Wang group of rank $n$}.
\item[(3)]  It is also well known (due to Mostow  \cite{MS}) that two solvmanifolds having
isomorphic fundamental groups are diffeomorphic.
\end{list}

It is often useful to assume that $k = b_1$ (the first Betti number of $M$) in the 
group extension of (1), which is possible due to a result  of Auslander and
Szczarba \cite{AS2} that a solvmanifold $M = D \backslash G$ has
{\em the canonical torus fibration} over the torus $ND \backslash G$ of dimension $b_1$
with fiber a nilmanifold, where $N$ is the nilradical of $G$
(the maximal connected normal subgroup of $G$).

It should be noted that a solvmanifold $M$ is not necessarily of the
form $\Gamma \backslash G$, where $G$ is a simply connected solvable Lie group
with discrete subgroup $\Gamma$. However, it is known (due to Auslander \cite{AU1})
that $M$ has a solvmanifold of the form $\Gamma' \backslash G$ as a finite covering,
where $\Gamma'$ is a subgroup of $\Gamma$ with finite index.
We know in general that two Wang groups are commensurable if and only if the
corresponding two solvmanifolds has the same solvmanifold as a finite covering.

\section{Four-dimensional solvmanifolds with complex structures}

Let $\Gamma$ be the fundamental group of a four-dimensional solvmanifold $S$.
Then we have
$$0 \rightarrow \Lambda \rightarrow \Gamma \rightarrow \Z^k \rightarrow 0,$$
where $\Lambda$ is a torsion-free nilpotent group of rank $4-k$.

For the classification of four-dimensional solvmanifolds (up to finite covering),
it is sufficient to classify Wang groups $\Gamma$ as the group extensions of
the above form, and find a subgroup $\Gamma'$ of finite index, which extends to a
simply connected solvable Lie group $G$ such that $\Gamma' \backslash G$ is a
solvmanifold.

We define three types of Wang groups of rank $4$ as follows:

\begin{list}{}{\topsep=5pt \leftmargin=25pt \itemindent=15pt}
\item[(Type I)] $2 \le k \le 4$,
\item[(Type II)] $k=1$ and $\Lambda$ is abelian,
\item[(Type III)] $k=1$ and $\Lambda$ is non-abelian,
\end{list}
where these three types are not mutually exclusive (as seen in Example 1).
\medskip

\noindent{\bfseries Example 1.}
Let $\Lambda_n= \Z^2 \rtimes \Z$, where the action $\phi:
\Z \rightarrow {\rm Aut}(\Z^2)$ is defined by
$\phi(1)=A_n \in {\rm GL}(2,\Z)$,
$$A_n = \left(
\begin{array}{cc}
1 & n\\
0 & 1
\end{array}
\right).$$
Then $\Lambda_n$ is a nilpotent group of rank $3$ and has a matrix expression:
$$\Lambda_n = \Bigg\{ \left(
\begin{array}[c]{ccc}
1 & a & \frac{c}{n}\\
0 & 1 & b\\
0 & 0 & 1
\end{array}
\right) \rule[-8mm]{0.25mm}{18mm}\;  a, b, c \in \Z \Bigg\}.
$$
We see that $\Lambda_n$ can be also expressed as a non-split group extension:
$$0 \rightarrow \Z \rightarrow \Lambda_n \rightarrow \Z^2 \rightarrow 0,$$
where the action of $\Z^2$ on $\Z$ is trivial. It should be noted that a
torsion-free nilpotent group of rank 3 is isomorphic to $\Lambda_n$ for
some $n \in \Z$. 

Let $\Gamma_n = \Lambda_n \times \Z$ be a nilpotent group of rank 4.
Then $\Gamma_n$ belongs to all of the types I, II and III.
\medskip

A four-dimensional solvmanifold $S$ is {\em of type I} ({\em II} or {\em III})
if the fundamental group of $S$ is of type I (II or III respectively).
Concerning the four-dimensional solvmanifolds of type I, we have
\medskip

\noindent {\bfseries Proposition 1.}  {\em A four-dimensional solvmanifold $S$ is
of type I if and only if $S$ is a $T^2$ bundle over $T^2$.}

{\em Proof.} We can easily see that Wang groups with $k=3$ or $4$ can be expressed
also as those with $k=2$. Hence the solvmanifolds of type I are all $T^2$ bundles
over $T^2$. Conversely, the fundamental group of a $T^2$ bundle over $T^2$ is
clearly a Wang group. Since we know \cite{SF} that diffeomorphism types of $T^2$
bundles over $T^2$ are also determined uniquely by their fundamental groups,
they are diffeomorphic to some solvmanifolds of dimension four. 
\enskip$\Box$
\medskip

It is known (due to Ue \cite{U}) that a complex
surface $S$ is diffeomorphic to a $T^2$ bundle over $T^2$ if and only if $S$ is a
complex torus,  Kodaira surface or hyperelliptic surface. The following result may be
considered as a generalization of this result.
\medskip

\noindent {\bfseries Theorem 1.}  {\em A complex surface is diffeomorphic to a
four-dimensional solvmanifold if and only if it is one of the following surfaces:
Complex torus, Hyperelliptic surface, Inoue Surface of type $S^0$,
Primary Kodaira surface, Secondary Kodaira surface, Inoue Surface of type $S^{\pm}$.
And every complex structure on each of these complex surfaces
(considered as solvmanifolds) is left-invariant.}
\medskip

The proof of Theorem 1 consists of three parts. The first part is to show that each of
the complex surfaces (in the theorem) can be characterized as a solvmanifold of type $II$
or $III$ with canonical complex structure. The rest of this section is devoted to this part of
the proof. The second part of the proof (the converse of the first part) is to show that a
complex surface with diffeomorphism type of solvmanifold must be one of the
complex surfaces in the theorem. We shall see this part of the proof in Section 6.
For the proof of the last part of the theorem, we shall give in Section 5 the complete
list of (left-invariant) complex structures on these complex surfaces as integrable
almost complex structures on their corresponding solvable Lie algebras.
\smallskip

We now characterize each of the complex surfaces in the theorem as a solvable
manifold of type $II$ or $III$ as follows:
\medskip

\noindent {[Type II]}\; We have the following split group extension, where
the action $\phi: \Z \rightarrow {\rm Aut} \Z^3$ is defined by
$\phi(1) \in {\rm SL}(3, \Z)$:
$$0 \rightarrow \Z^3 \rightarrow \Gamma \rightarrow \Z \rightarrow 0.$$

\begin{list}{}{\topsep=5pt \leftmargin=0pt \itemindent=20pt \itemsep=3pt}
\item[(1)] {\em Complex Tori}
\begin{quote} $\phi(1)=I$, and thus $\Gamma=\Z^4$
\end{quote}
\item[(2)] {\em Hyperelliptic Surfaces}
\begin{quote} $\phi(1)$ has a single root $1$, and a double root $-1$ with
linearly independent eigenvectors or non-real complex roots
$\beta, \bar{\beta}$ with $|\beta|=1$.
\end{quote}
\item[(3)] {\em Inoue Surfaces of type $S^0$}
\begin{quote} $\phi(1)$ has non-real complex roots $\alpha, \bar{\alpha}$
and  a real root $c\;(c \not=1)$ with $|\alpha|^2c = 1$.
\end{quote}
\end{list}
\medskip

\noindent [Type III] \; The group extension is split, and thus determined
only by the action $\phi: \Z \rightarrow {\rm Aut}(\Lambda_n)$:
$$0 \rightarrow \Lambda_n \rightarrow \Gamma \rightarrow \Z \rightarrow 0,$$
where $\Lambda_n$ is a nilpotent group of rank $3$ as defined in Example 1.

An automorphism $\phi(1)$ induces the automorphism $\widetilde{\phi}(1)$ of the
center $\Z$ of $\Lambda_n$, and the automorphism $\widehat{\phi}(1)$ of
$\Z^2 = \Lambda_n/\Z$.
\medskip

\begin{list}{}{\topsep=5pt \leftmargin=0pt \itemindent=20pt \itemsep=3pt}
\item[(4)] {\em Primary Kodaira Surfaces}
\begin{quote} $\phi(1) = {\rm Id}$, and thus $\Gamma=\Lambda_n \times \Z$.
\end{quote} 
\item[(5)] {\em Secondary Kodaira Surfaces}
\begin{quote} $\widetilde{\phi}(1)= {\rm Id}$, and $\widehat{\phi}(1)$ has a double root $1$
with linearly independent eigenvectors or non-real complex roots
$\alpha, \bar{\alpha}$ with $|\alpha|=1$.
\end{quote}
\item[(6)] {\em  Inoue Surfaces of Type $S^{\pm}$}
\begin{quote} $\widetilde{\phi}(1)= \pm {\rm Id}$, and $\widehat{\phi}(1)$ has two positive
roots $a, b\; ( a \not=1))$ with $ab=1$ (two real roots $a, b$ with opposite sign
$(|a| \not=1)$ with $ab=-1$ respectively).
Note that the Inoue surface of type $S^-$ has that of type $S^+$ as a double covering.
\end{quote}
\end{list}
\medskip

\noindent (1) {\bfseries \emph{Complex Tori}}
\smallskip

An $n$-dimensional torus $T^n$ is a compact homogeneous space of the abelian
Lie group $\R^n$: that is, $T^n = \Z^n \backslash \R^n$ where $\Z^n$ is an abelian
lattice of $\R^n$ which is spanned by some basis of $\R^n$ as a real vector space.
For the case $n = 2m$, the standard complex structure $\C^m$ on $\R^{2m}$ defines
a complex structure on $T^{2m}$. The complex manifold thus obtained is a
{\em complex torus}.
\smallskip

It should be noted that $T^n (n \ge 3)$ can admits a
structure of non-toral solvmanifold (see Section 7).
On the other hand it is unknown if all the complex
structures on $T^{2m}\; (m \ge 3)$ are the standard ones.
This holds for $m = 1, 2$, since
it does define a Riemann surface with $b_1=2$ (a complex surface with $b_1= 4$
respectively), which is K\"ahlerian, and its albanese map is biholomophic.
\medskip

\noindent (2) {\bfseries \emph{Hyperelliptic Surfaces}}
\smallskip

Let $\Gamma = \Z^3 \rtimes \Z^1$, where the action
$\phi: \Z^1 \rightarrow {\rm Aut}(\Z^3)$ is defined by
$\phi(1) \in {\rm SL}(3,\Z)$. Assume that $\phi(1)$ has a single root $1$,
and a double root $-1$ with linearly independent eigenvectors of $-1$
or non-real complex roots $\beta, \bar{\beta}$ with $|\beta|=1$.
\smallskip

For $A = \phi(1) \in {\rm SL}(3,\Z)$ which satisfies our assumption,
we can find a basis $\{u_1, u_2, u_3\}$ of $\R^3$
such that $A u_1 = au_1 - bu_2, A u_2 = bu_1 + au_2, A u_3 = u_3$,
where $a = -1$ and $b = 0$ for the case that $A$ has a double root $-1$,
and $a ={\rm Re}\,\beta$ and $b = {\rm Im}\,\beta$ for the case that $A$ has a
non-real complex root $\beta$. 
Let $u_i=(u_{i1}, u_{i2}, u_{i3}),\; i=1,2,3$.
Then  $\{v_1, v_2, v_3 \},\;v_j=(u_{1j}, u_{2j}, u_{3j}), j=1,2,3$ defines an
abelian lattice $\Z^3$ of $\R^3$ which is preserved by a rotation around a fixed
axis. In particular, $\beta$ must be $e^{\sqrt{-1} \zeta}$
($\zeta = \frac{2}{3} \pi, \frac{1}{2} \pi$ or $\frac{1}{3} \pi$).

Furthermore, we may assume that $u_{3j} = 0, j = 1,2,$ and $A$ is of the form
$$\left(
\begin{array}[c]{ccc}
a_{11} & a_{12} & 0\\
a_{21} & a_{22} & 0\\
p & q & 1
\end{array}
\right),$$
where $A' = (a_{ij}) \in {\rm SL}(2, \Z), p, q \in \Z$.
Since $A'$ has  the root $-1$ (with linearly independent eigenvectors) or $\beta$,
we can assume that $A'$ is of the form: 
$$\left(
\begin{array}[c]{cc}
-1 & 0\\
0 & -1
\end{array}
\right),
\left(
\begin{array}[c]{cc}
0 & 1\\
-1 & -1
\end{array}
\right),
\left(
\begin{array}[c]{cc}
0 & 1\\
-1 & 0
\end{array}
\right),
\left(
\begin{array}[c]{cc}
0 & 1\\
-1 & 1
\end{array}
\right),
$$
according to the root $e^{\sqrt{-1}\, \eta}$ of $A'$, where
$\eta = \pi, \frac{2}{3} \pi, \frac{1}{2} \pi$ or $\frac{1}{3} \pi$, respectively. 
\smallskip

We now define a solvable Lie group $G = (\C \times \R) \rtimes \R$, where the
action $\phi: \R \rightarrow {\rm Aut}(\C \times \R)$ is defined by 
$$\phi(t)((z,s)) = (e^{\sqrt{-1} \eta t} z, s),$$
which is a canonical extension of $\phi$.
$\Gamma = \Z^3 \rtimes \Z$ clearly defines a lattice of $G$. 
Since the action on the second factor $\R$ is trivial, the multiplication of $G$ is
defined on $\C^2$ as follows:
$$(w_1,w_2) \cdot (z_1,z_2)= (w_1+e^{\sqrt{-1} \eta t} z_1, w_2+z_2),$$
where $t= {\rm Re}\,w_2$.  For each lattice $\Gamma$ of $G$,
$S=\Gamma \backslash G$ with the canonical complex structure from $\C^2$
defines a complex surface which is, by definition, a {\em hyperelliptic surface} 
(see \cite{BPV}). 

We can see that there exist seven
isomorphism classes of lattices of $G$, which correspond to seven classes of
hyperelliptic surfaces.
For each $\eta$, take a lattice $\Z^3$ of $\R^3$ spanned by $\{v_1,v_2,v_3\}$ for
$A$ with $p = q = 0$. Then we can get a lattice $\Z^3$ spanned by $\{v_1,v_2,v_3'\}$
for $A$ with arbitrary $(p,q) \in \Z^2$, by changing $v_3$ into
$v_3' = s v_1 + t v_2 + v_3$ where $s, t \in \Q$ with $0 \le s,t < 1$, and
$\Gamma = \Z^3 \rtimes \Z$ defines a lattice of the solvable Lie group $G$.
By elementary calculation, we obtain the following seven isomorphism classes of
lattices: besides four trivial cases with $(p,q) = (0,0)$ and $(s,t) = (0,0)$ for
$\eta = \pi, \frac{2}{3} \pi, \frac{1}{2} \pi$ and $\frac{1}{3} \pi$,
we have three other cases with $(p,q) = (1,0)$, and 
(i) $(s,t) = (\frac{1}{2}, 0)$ for $\eta = \pi$,
(ii) $(s,t) = (\frac{1}{3}, \frac{1}{3})$ for $\eta = \frac{2}{3} \pi$,
(iii) $(s,t) = (\frac{1}{2}, \frac{1}{2})$ for $\eta = \frac{1}{2} \pi$.
\medskip

\noindent (3) {\bfseries \emph{Inoue Surfaces of Type ${\bf S}^0$}}
\smallskip

Let $\Gamma= \Z^3 \rtimes \Z$, where the action $\phi:
\Z \rightarrow {\rm Aut}(\Z^3)$ is defined by $\phi(1) \in {\rm SL}(3,\Z)$.
Assume that $\phi(1)$ has complex roots $\alpha, \bar{\alpha}$ and a real root
$c\; (c \not= 1)$ with $|\alpha|^2 c= 1$.

Let $(\alpha_1,\alpha_2,\alpha_3) \in \C^3$ be the eigenvector of $\alpha$ and
 $(c_1,c_2,c_3) \in \R^3$ the eigenvector of $c$. The set of vectors
$\{(\alpha_i,c_i) \in \C \times \R \,|\, i= 1,2,3\}$ are linearly independent over $\R$,
and defines a lattice $\Z^3$ of $\C \times \R$. 
Then $\Gamma=\Z^3 \rtimes \Z$ can be extended to a solvable Lie group
$G= (\C \times \R) \rtimes \R$,
where the action $\bar{\phi}: \R \rightarrow {\rm Aut}(\C \times \R)$ is defined by
$\bar{\phi}(t): (z,s) \rightarrow (\alpha^t z,c^t s)$.
\medskip

Taking a coordinate change $\R \rightarrow \R_+$ defined by $t \rightarrow
e^{\log c t}$ and regarding $\R \times \R_+$ as ${\bf H}$ (the upper half plain),
$M=\Gamma \backslash G$ can be considered as
$\C \times {\bf H}/\Gamma'$, where $\Gamma'$ is
a group of automorphisms generated by $g_0$ and $g_i, i=1,2,3$,
which correspond to the canonical generators of $\Gamma$.
To be more precise, we see
$$g_0:(z_1,z_2) \rightarrow (\alpha z_1, c z_2),
g_i:(z_1,z_2) \rightarrow (z_1 + \alpha_i, z_2 + c_i), i=1,2,3.$$
$S=\C \times {\bf H}/\Gamma'$ is, by definition, an {\em Inoue surface of type $S^0$}
\cite{I}.
\medskip

\noindent (4) {\bfseries \emph{Primary Kodaira Surfaces}}
\smallskip

Let $\Gamma_n = \Lambda_n \times \Z$ (a nilpotent group of rank 4).
$\Gamma_n$ can be extended to the nilpotent Lie group $G = N \times \R$,
where 
$$N = \Bigg\{ \left(
\begin{array}[c]{ccc}
1 & x & s\\
0 & 1 & y\\
0 & 0 & 1
\end{array}
\right) \rule[-8mm]{0.25mm}{18mm}\;  x, y, s \in \R \Bigg\}.
$$
Taking the coordinate change $\Phi$ from $N \times \R$ to $\R^4$:
$$\Phi: ((x,y,s), t) \longrightarrow
(x,y, 2 s-x y, 2 t+ \frac{1}{2}(x^2+y^2)),$$
and regarding $\R^4$ as $\C^2$, the group operation on $G$ can be expressed as
$$(w_1,w_2) \cdot (z_1,z_2)=
(w_1 + z_1, w_2 - \sqrt{-1} \bar{w_1} z_1 + z_2).$$
 Let $\Gamma_n'$ is the corresponding group of affine transformations on $\C^2$,
then $S= \C^2/\Gamma_n'$ is, by definition, a {\em Primary Kodaira surface} \cite{K1}.
\medskip

\noindent (5) {\bfseries \emph{Secondary Kodaira surfaces}}
\smallskip

Let $\Gamma_n= \Lambda_n \rtimes \Z$, where the action $\phi:
\Z \rightarrow {\rm Aut}(\Lambda_n)$ satisfies the condition that
the induced automorphism $\widetilde{\phi}(1)$ of $\Z$ is trivial, that is,
$\widetilde{\phi}(1)= {\rm Id}$, and the induced automorphism
$\widehat{\phi}(1)$ of $\Z^2$ has a double root $-1$
with linearly independent eigenvectors, or non-real complex roots,
$\alpha, \bar{\alpha}\,
(\alpha \not= \bar{\alpha})$ with $|\alpha|= 1$.

We shall see that $\Gamma_n$ can be extended to a solvable Lie group
$G = N \rtimes \R$. As we have seen in (2), $\alpha$ must be
$e^{i \eta}$, $\eta = \pi, \frac{2}{3} \pi, \frac{1}{2} \pi$ or $\frac{1}{3} \pi$,
and there exists a basis $\{u_1', u_2'\}$ of $\R^2$ such that
$A u_1' = a u_1' - b u_2', A u_2' = b u_1' + a u_2'$, where $A = \widehat{\phi}(1)$,
$a = {\rm Re}\,\alpha, b = {\rm Im}\,\alpha$,
and $u_1' = (u_{11}, u_{12}), u_2' = (u_{21}, u_{22})$. 
The abelian lattice $\Z^2$ of $\R^2$ spanned by $\{v_1', v_2'\}$,
where $v_1' = (u_{11}, u_{21}), v_2' = (u_{12}, u_{22})$, is preserved by
the automorphism $\psi': (x, y) \rightarrow (ax - by, bx + ay)$ of $\R^2$.
We can extend $\psi'$ to an automorphism of $N$ of the form:
$$\psi: (x,y,\frac{z}{n}) \longrightarrow (a x-b y, b x+ a y, \frac{z}{n}+ h(x, y)),$$
where  $h(x, y) = \frac{1}{2} b\,(a x^2 - a y^2 - 2 b x y)$.
We can extend the lattice $\Z^2$ spanned by $\{v_1', v_2'\}$ to a lattice
$\Lambda_n$ spanned by $\{v_1, v_2, v_3\}, v_1 = (u_{11}, u_{21}, u_{31}),
v_2 = (u_{12}, u_{22}, u_{32}), v_3 = (0, 0, u_{33})$ for suitable $u_{31}, u_{32},
u_{33}$, so that $\Lambda_n$ is preserved by $\psi$. We now define a solvable
Lie group $G$ by extending the action $\phi(m) = \psi^m, m \in \Z$ to
$\psi(t), t \in \R$, replacing $a$ with $\cos \eta t$ and $b$ with $\sin \eta t$.
It is clear that $\Gamma_n = \Lambda_n \rtimes \Z$ defines a lattice of $G$.
If we take the new coordinate as in (4), the automorphism $\psi$ is
expressed as,
$$(z_1, z_2) \longrightarrow (\zeta z_1, z_2)$$
for $\zeta \in \C, |\zeta|=1$. It follows that the above automorphism
is holomorphic with respect to the complex structure defined in (3).
We see that $S_n= \Gamma_n \backslash G$ is a
finite quotient of a primary Kodaira surface; and $S_n$
with the above complex structure is, by definition,
{\em a secondary Kodaira surface} \cite{BPV, FM}.
\medskip

\noindent (6) {\bfseries \emph{Inoue surfaces of type $\bf S^\pm$}}
\smallskip

Let $\Gamma_n= \Lambda_n \rtimes \Z$, where the action $\phi:
\Z \rightarrow {\rm Aut}(\Lambda_n)$ satisfies the condition that
for the induced action $\widetilde{\phi}: \Z \rightarrow {\rm Aut}(\Z)$,
$\widetilde{\phi}(1)= {\rm Id}$, and for the induced action
$\widehat{\phi}: \Z \rightarrow {\rm Aut}(\Z^2)$,
$\widehat{\phi}(1)= (n_{ij}) \in {\rm SL}(2, \Z)$ has two positive real roots
$a, b$ with $ab=1$.

Let $(a_1,a_2), (b_1,b_2) \in \R^2$ be
eigenvectors of $a, b$ respectively. Let $G= N \rtimes \R$ be a solvable
Lie group, where the action $\bar{\phi}: \R \rightarrow {\rm Aut}(N)$
is defined by
$$\bar{\phi}(t): \left(
\begin{array}[c]{ccc}
1 & x & z\\
0 & 1 & y\\
0 & 0 & 1
\end{array}
\right)
\longrightarrow
\left(
\begin{array}[c]{ccc}
1 & a^t x & z\\
0 & 1 & b^t y\\
0 & 0 & 1
\end{array}
\right),
$$
which is a canonical extension of $\phi$.
In order to define a lattice $\Lambda_n$ which is preserved by
$\bar{\phi}$, we take $g_1,g_2,g_3 \in N$ as
$$g_1=\left(
\begin{array}[c]{ccc}
1 & a_1 & c_1\\
0 & 1 & b_1\\
0 & 0 & 1
\end{array}
\right),\;
g_2=\left(
\begin{array}[c]{ccc}
1 & a_2 & c_2\\
0 & 1 & b_2\\
0 & 0 & 1
\end{array}
\right),\;
g_3=\left(
\begin{array}[c]{ccc}
1 & 0 & c_3\\
0 & 1 & 0\\
0 & 0 & 1
\end{array}
\right),$$
where $c_1,c_2,c_3$ are to be determined, satisfying
the following conditions:

\begin{list}{}{\partopsep=0pt \parsep=0pt \leftmargin=20pt}
\item[1)]  $[g_1, g_2]=g_3^n$
\item[2)]  $\bar{\phi}(1)(g_1)=g_1^{n_{11}} g_2^{n_{12}} g_3^k,\;
\bar{\phi}(1)(g_2)=g_1^{n_{21}} g_2^{n_{22}} g_3^l$, where $k,l \in \Z.$
\end{list}
If we take $g_0 \in N \rtimes \R$ as
$$g_0=\left( \left(
\begin{array}[c]{ccc}
1 & 0 & p\\
0 & 1 & 0\\
0 & 0 & 1
\end{array}
\right), 1
\right),
\;p \in \R,$$
then $\{g_0,g_1,g_2,g_3\}$ defines a lattice $\Gamma_n$ of $G$, and
$S_n=\Gamma_n\backslash G$ is a solvmanifold.
\smallskip

Now, we define a diffeomorphism $\Psi: G=N \rtimes
\R \longrightarrow \R^3 \times \R_+$, for an arbitrary
$\gamma=p + q \sqrt{-1} \in \C$ and $\sigma=\log b$, by
$$\Psi: \left( \left(
\begin{array}[c]{ccc}
1 & y & x\\
0 & 1 & s\\
0 & 0 & 1
\end{array}
\right),
t \right)
\longrightarrow
(x,\, e^{\sigma t} y+ q \,t,\, s, \,e^{\sigma t}).$$
Then considering $\R^3 \times \R_+$ as $\C \times {\bf H}$,
$g_0,g_1,g_2,g_3$ are corresponding to the following holomorphic
automorphisms of $\C \times {\bf H}$,
$$g_0: (z_1, z_2) \rightarrow (z_1 + \gamma, b z_2),$$
$$g_i: (z_1, z_2) \rightarrow (z_1 + a_i z_2 + c_i, z_2 + b_i),$$
where $i=1,2, 3$ and $a_3 = b_3 = 0$.
$S_n$ with the above complex structure is, by definition,
an {\em Inoue surface of type $S^+$} \cite{I}.
\smallskip

An {\em Inoue surface of type $S^-$} is defined similarly as
the case where the action $\phi: \Z \rightarrow {\rm Aut}(\Lambda_n)$
satisfies the condition that $\widetilde{\phi}(1)= -{\rm Id}$,
and $\widehat{\phi}(1)$ has a positive and a negative real root.
It is clear that an Inoue surface of type $S^-$ has $S^+$ with $\gamma=0$
as its double covering surface.

\section{Examples}

We give in this section three other classes of four-dimensional (orientable)
solvmanifolds which (as a consequence of Theorem 1)
admit no complex structures.
\medskip

\noindent{\bfseries Example 2.}  Let $\Gamma$ be the Wang group of type II defined
by the following split group extension, where
the action $\phi: \Z \rightarrow {\rm Aut} \Z^3$ is defined by
$\phi(1) \in {\rm SL}(3, \Z)$:
$$0 \rightarrow \Z^3 \rightarrow \Gamma \rightarrow \Z \rightarrow 0.$$

\begin{list}{}{\topsep=5pt \leftmargin=10pt \itemindent=10pt \itemsep=3pt}
\item[(1)] Suppose that $\phi(1)$ has three distinct positive real roots
$a_1,a_2,a_3$, then there exist linearly independent eigenvectors
$u_1, u_2, u_3$ of $a_1, a_2, a_3$ respectively.
Let $u_i=(u_{i1}, u_{i2}, u_{i3}),\; i=1,2,3$. Then we have an abelian lattice
$\Z^3$ of $\R^3$ defined by
$\{v_1, v_2, v_3 \},\;v_j=(u_{1j}, u_{2j}, u_{3j}), j=1,2,3$.
We define a solvable Lie group
$G = \R^3 \rtimes \R$,
where the action $\phi: \R \rightarrow \Aut(\R^3)$ is defined by
$$\phi(t)((x, y, z)) = (e^{t \log a_1} x, e^{t \log a_2} y, e^{t \log a_3} z),$$
which is a canonical extension of $\phi$.
Then $\Gamma = \Z^3 \rtimes \Z$ is a lattice of $G$,
and  $S=\Gamma \backslash G$ is a solvmanifold. We can see that
$S$ is a $T^2$-bundle over $T^2$ with $b_1=2$ for the case
where one of the roots is $1$, and a $T^3$ bundle over $T^1$ with $b_1=1$ for
the case where none of the roots is $1$.

\item[(2)] Suppose that  $\phi(1)$ has a triple root $1$, then taking a suitable basis
$\{u_1,u_2,u_3\}$ of $\R^3$, $\phi(1)$ is expressed in either of the following
forms:
$$\left(
\begin{array}[c]{ccc}
1 & 1 & \frac{1}{2}\\
0 & 1 & 1\\
0 & 0 & 1
\end{array}
\right),\;
\left(
\begin{array}[c]{ccc}
1 & 1 & 0\\
0 & 1 & 0\\
0 & 0 & 1
\end{array}
\right).
$$
Let $G=\R^3 \rtimes \R$, where the action $\bar{\phi}: \R \rightarrow
{\rm Aut}(\R^3)$ is defined by
$$\bar{\phi}(t)=
{\rm exp}\;t
\left(
\begin{array}[c]{ccc}
0 & 1 & 0\\
0 & 0 & 1\\
0 & 0 & 0
\end{array}
\right)=\left(
\begin{array}[c]{ccc}
1 & t & \frac{1}{2} t^2\\
0 & 1 & t\\
0 & 0 & 1
\end{array}
\right)
$$
for the former case, and
$$\bar{\phi}(t)=
{\rm exp}\;t
\left(
\begin{array}[c]{ccc}
0 & 1 & 0\\
0 & 0 & 0\\
0 & 0 & 0
\end{array}
\right)=\left(
\begin{array}[c]{ccc}
1 & t & 0\\
0 & 1 & 0\\
0 & 0 & 1
\end{array}
\right)
$$
for the latter case.
Then, as defined in the above case (1), $\{v_0, v_1, v_2, v_3\}$ defines a
lattice $\Gamma$ of $G$, and $S=\Gamma \backslash G$ is a nilmanifold.
We can see that $S$ is a nilmanifold with $b_1=2$ for the former case (which
admits no complex structures), and
a nilmanifold with $b_1=3$ for the latter case (which is a primary Kodaira surface).
\end{list}
\medskip

\noindent {\bfseries Remark.} We can show that, up to finite covering,
there exit nine classes of four-dimensional (orientable) solvmanifolds:
six classes of complex surfaces in Section 3, and three classes of solvmanifolds
in Example 2 which admit no complex structures.
\medskip

\noindent{\bfseries Example 3.}  Let $\Gamma_n$ be a Wang group expressed as the
extension:
$$0 \rightarrow \Z^2 \rightarrow \Gamma_n \rightarrow \Z^2 \rightarrow 0,$$

\noindent where the action $\phi: \Z^2 \rightarrow \Aut(\Z^2)$ is defined by
$\phi(e_1), \phi(e_2), e_1 = (1,0), e_2 = (0,1)$. Suppose that
$\phi(e_1) = -{\rm I}$, and $\phi(e_2)$ is of the form
$$\left(
\begin{array}{cc}
1 & n\\
0 & 1
\end{array}
\right).$$
Then we can express $\Gamma_n$ also as $\Lambda_n \rtimes \Z$, where $\Lambda_n$
is the nilpotent group of rank $3$ as defined in Example 1, and the action $\phi:
\Z \rightarrow {\rm Aut}(\Lambda_n)$ is defined by $\phi(1)= \tau
\in {\rm Aut}(\Lambda_n)$,
$$\tau:  \left(
\begin{array}[c]{ccc}
1 & x & \frac{z}{n}\\
0 & 1 & y\\
0 & 0 & 1
\end{array}
\right) \;
\longrightarrow \;
\left(
\begin{array}[c]{ccc}
1 & x & \frac{-z}{n}\\
0 & 1 & -y\\
0 & 0 & 1
\end{array}
\right)
$$
\smallskip

Let $G= N \times \R$ be a nilpotent Lie group, where $N$ is the
nilpotent Lie group obtained by the real completion of $\Lambda_n$ as defined in (4).
$\Gamma_n$  acts freely as a group of automorphisms on $G$, and
$S_n=G/\Gamma_n$ is a solvmanifold with $b_1=2$.
$S_n$ has the  nilmanifold $(\Lambda_n \times \Z) \backslash G$ as a double
covering, with the covering transformation group $\Z_2$ generated by $\tau$.
\smallskip

We show that $S_n$ is parallelizable, that is, $S_n$ admits a field of linear frame
(consisting of four linearly independent vector fields).
We have three linearly independent left-invariant vector fields on $N$:
$$X_1=\frac{\partial}{\partial x}, X_2= \frac{\partial}{\partial y}+nx \frac{\partial}{\partial z},
X_3= \frac{\partial}{\partial z}.$$ 
We define vector fields on $G$ invariant by
$\tau$: $\widetilde{X}_1=X_1, \widetilde{X}_4=\frac{\partial}{\partial t}$, and 
$$\widetilde{X}_2=\cos (\pi t) X_2 + \sin (\pi t) X_3,
\widetilde{X}_3=-\sin (\pi t) X_2 + \cos (\pi t) X_3, $$
which are linearly independent and invariant by $\Gamma_n$.
\smallskip

It should be noted that $S_n$ can be expressed as $D_n \backslash \overline{G}$,
where $\overline{G}$ is the extension of $G$ with $y, z \in \C$ and the action
$\overline{\phi}$ defined by
$$\overline{\phi}(t) (x, y, z)= (x, e^{\sqrt{-1} \pi t} y, e^{\sqrt{-1} \pi t} z),$$
and $D_n=\Gamma_n H$,  where $H$ is the closed subgroup of $\overline{G}$ with
$x = 0$, $y = \sqrt{-1} y_2$, $z = \sqrt{-1} z_2$, $y_2, z_2 \in \R$.
\medskip

\noindent {\bfseries Remark.} As noted in the paper \cite{AS1}, all of the
four-dimensional solvmanifolds are
parallelizable. This is trivial for the case where $S$ is of the form $\Gamma \backslash G$
(where $\Gamma$ is a simply connected solvable Lie group with lattice $\Gamma$).
For the case where $S$ is not of the form  $\Gamma \backslash G$, we can see that, as
shown in  Example 3, $S$ always admit a field of linear frame.

\section{ Complex structures on solvable Lie algebras}

A four-dimensional solvmanifold can be written
(up to finite covering) as $\Gamma \backslash G$, where $\Gamma$ is a lattice
of a simply connected solvable Lie group $G$.  In this section we express complex
structures of the complex surfaces in the last section as those induced from
left-invariant complex structures on $G$ (i.e. left-invariant complex
structures), by defining integrable almost complex structures $J$ on the Lie algebra
 $\mathfrak g$ of $G$.
In the following list, we express the Lie algebra $\mathfrak g$ of $G$ as having
a basis $\{X_1,X_2,X_3,X_4\}$  with the bracket multiplication specified for each
complex surface.
Except for (6), the almost complex structure $J$ is defined by
$$JX_1=X_2, JX_2=-X_1, JX_3=X_4, JX_4=-X_3,$$
for which the Nijenhuis tensor
$N_J(X_i,X_j)=[JX_i,JX_j]-J[JX_i,X_j]-J[X_i,JX_j]-[X_i,X_j], 1 \le i < j \le 4$ vanishes.
\medskip

{\baselineskip=16pt
\begin{list}{}{\topsep=0pt \leftmargin=5pt \itemindent=15pt \parsep=0pt \itemsep=3pt}
\item[ (1)] Complex Tori \par
$[X_i, X_j]=0 \;( 1 \le i < j \le 4)$.
\item[ (2)] Hyperelliptic Surfaces \par
$[X_4, X_1]= -X_2, [X_4, X_2]= X_1$, and all other brackets vanish.
\item[ (3)] Primary Kodaira Surfaces \par
$[X_1, X_2]= -X_3$, and all other brackets vanish.
\item[ (4)] Secondary Kodaira Surfaces \par
$[X_1, X_2]= -X_3, [X_4, X_1]= -X_2, [X_4, X_2]= X_1$,
and all other brackets vanish.
\item[ (5)] Inoue Surfaces of Type $S^0$ \par
$[X_4, X_1]= a X_1 - b X_2, [X_4, X_2]= b X_1 + a X_2, [X_4, X_3]= -2a X_3$,
and all other brackets vanish, where $a, b \;(\not=0) \in \R$.
\item[ (6)] Inoue Surfaces of Type $S^+$ and $S^-$ \par
$[X_2, X_3]= -X_1, [X_4, X_2]= X_2, [X_4, X_3]= -X_3$, and all other
brackets vanish.
The almost complex structure $J$ is defined by
$$J X_1=X_2, J X_2=-X_1, J X_3= X_4-q X_2, J X_4= -X_3-q X_1,$$
for which the Nijenhuis tensor $N_J$ vanishes.
\end{list}
}

\section{Proof of Theorem 1}

We have seen in Section 3 and 5 that each of the complex surfaces in Theorem 1
can be characterized as a four-dimensional solvmanifold with canonical
left-invariant complex structure. In this section we shall complete the proof of Theorem 1,
showing the converse that a complex surface with diffeomorphism of solvmanifold
must be one of the complex surfaces in the theorem.

We denote by $S$ a complex surface with diffeomorphism type of solvmanifold.
We first remark that since $S$ is parallelizable,
the Euler number $c_2$ of $S$ vanishes, and the fundamental group of $S$
is abelian if and only if $S$ is a four-dimensional torus \cite{MS}.
Let $\kappa(S)$ be the Kodaira dimension of $S$. The classification of complex
surfaces with $c_2 = 0$ is divided into three cases: $\kappa(S) = -\infty, 0, 1$
(\cite{BPV}). In the case where $\kappa(S) = -\infty$, $S$ is a surface of class
${\rm VII}_0$ or ruled surface of genus 1. The latter surface cannot be diffeomorphic
to a solvmanifold since the fundamental group of ruled surface of genus $1$
is $\Z^2$ (which is abelian).
According to the well-known theorem of Bogomolov 
(proved by Li, Yau and Zheng \cite{LYZ}), we know that a complex surface of class
${\rm VII_0}$ with $b_1 = 1$ and $c_2 = 0$ is an Inoue surface or a Hopf surface. 
Since the fundamental group of the Hopf surface is of the form $H \rtimes \Z$
where $H$ is a finite unitary group (including the trivial case),
it cannot be the fundamental group of solvmanifold. 
Hence $S$ must be an Inoue surface. In the case where $\kappa(S) = 0$,
$S$ is a complex torus, hyperelliptic surface, or Kodaira surface (primary or
secondary Kodaira surface).

In the case where $\kappa(S) = 1$, $S$ is, a (properly) elliptic surface (which is
minimal since $c_2 = 0$). Let us first recall some terminologies and fundamental
results concerning topology of elliptic surfaces in general. An elliptic surface is a
complex surface $S$ together with an elliptic fibration $f: S \rightarrow B$ where
$B$ is a curve, such that a general fiber $f^{-1}(t), t \in B$
(except finite points $t_1, t_2, ..., t_k$) is an elliptic curve. The base curve $B$
is regarded as a two-dimensional orbifold with multiple points $t_i$ with
multiplicity $m_i$, where $m_i \,(i \ge 2)$ is the multiplicity of the fiber
$f^{-1}(t_i)$ ($i = 1, 2, ..., k$). An elliptic surface $S$ is of the type {\em hyperbolic,
flat (Euclidean), spherical or bad}, according as the orbifold $B$ is of that type.
The Euler number $e^{orb}(B)$ of $B$ is by definition
$e(B) - \sum_{i = 1}^k (1 - \frac{1}{ m_i})$, where $e(B)$ is the Euler number of $B$
as topological space. We know that $B$ is hyperbolic, flat or spherical according
as $e^{orb}(B)$ is negative, 0, or positive. In the case $c_2 = 0$, we can see
\cite{W1, W2} that an elliptic surface $S$ is hyperbolic, flat or spherical according
as the Kodaira dimension $\kappa(S)$ is $1, 0$ or $-\infty$. 
We now continue our proof for the case $\kappa(S) = 1$. By the above argument,
$S$ is a minimal elliptic surface of hyperbolic type. We show that the
fundamental group of $S$ is not solvable, and thus $S$ cannot be diffeomorphic to
a solvmanifold. We have the following presentation of $\pi_1(S)$ as a short
exact sequence \cite{FM}:
$$ 0 \rightarrow \Z^2 \rightarrow \pi_1(S) \rightarrow \pi_1^{orb}(B) \rightarrow 0,$$

\noindent where $\pi_1^{orb}(B)$ is the fundamental group as two-dimensional
orbifold.

Since $\pi_1^{orb}(B)$ is a discrete subgroup of $\rm PSL(2, \R)$, it contains
a torsion-free subgroup $\Gamma$ of finite index, such that $\Gamma$ is the
fundamental group of a finite orbifold covering $\widetilde{B}$ of $B$, which is a
closed surface of genus $g \ge 2$ (or a Riemann surface of hyperbolic type)
\cite{SC}. We know that $\Gamma$ is represented as a group with generator
$\{a_1, a_2, \ldots, a_g, b_1, b_2, \ldots, b_g\}$ and relation
$\prod_{i = 1}^g [a_i, b_i] = 1$, which is not solvable for $g \ge 2$ (\cite{CH}).
It follows that $\pi_1(S)$ cannot be solvable since the quotient groups and
subgroups of a solvable group must be solvable. This completes the proof of
Theorem 1. 
\hfill $\Box$

\section{Solvmanifolds with K\"ahler structures}

Let $M$ be a solvmanifold of the form $ \Gamma \backslash G$,
where $\Gamma$ is a lattice of a simply connected solvable Lie group $G$.
$M$ is {\em of completely solvable type},
if the adjoint representation of the Lie algebra $\mathfrak g$ of $G$ has
only real eigenvalues; and {\em of rigid type} (or {\em of type (R)})
in the sense of Auslander, if the adjoint representation of $\mathfrak g$ has only
pure imaginary (including $0$) eigenvalues.
It is clear that $M$ is both of completely solvable and of rigid type if and
only if $\mathfrak g$ is nilpotent, that is, $M$ is a nilmanifold.
As we have seen in Section 4 and 5, a hyperelliptic surface can be characterized as
a solvmanifold of rigid type with canonical complex structure. we define a natural
generalization of hyperelliptic surfaces as in the following example.
\medskip

\noindent {\bfseries Example 4.} 
Let $G = \C^l \rtimes \R^{2k}$ with the action
$\phi: \R^{2k} \rightarrow {\rm Aut}(\C^l)$ defined by
$$\phi(\bar{t}_i) ((z_1, z_2, \ldots, z_l)) = 
(e^{\sqrt{-1}\,\eta^i_1\, t_i} z_1, e^{\sqrt{-1}\,\eta^i_2\, t_i} z_2, ...,
e^{\sqrt{-1}\, \eta^i_l\, t_i} z_l),$$
where $\bar{t}_i= t_i e_i$ ($e_i$: the $i$-th unit vector in $\R^{2k})$, and
$e^{\sqrt{-1}\,\eta^i_j}$ is the $s_i$-th root of unity,
$i = 1, \ldots,2k, j = 1, \ldots,l$. If an abelian lattice $\Z^{2l}$ of $\C^l$ is preserved by
the action $\phi$ on $\Z^{2k}$, then $M=\Gamma \backslash G$ defines a solvmanifold
of rigid type, where $\Gamma=\Z^{2l} \rtimes \Z^{2k}$ is a lattice of $G$.
In fact the Lie algebra $\mathfrak g$ of $G$ is the following:
$$ \mathfrak g = \{X_1, X_2, \ldots , X_{2l}, X_{2l+1}, \ldots , X_{2l+2k}\},$$
where the bracket multiplications are defined by 
$$[X_{2l+2i}, X_{2j-1}] = -X_{2j}, [X_{2l+2i}, X_{2j}] = X_{2j-1}$$
for $i = 1, \ldots,k, j = 1, \ldots, l$,
and all other brackets vanish. The canonical left-invariant complex structure is defined by
$$JX_{2j-1}=X_{2j}, JX_{2j}=-X_{2j-1}, JX_{2l+2i-1}=X_{2l+2i}, JX_{2l+2i}=-X_{2l+2i-1}$$ 
for $i=1, \ldots, k, j=1, \ldots, l$.
\medskip

\noindent {\bfseries Remark.} Let $G = \C^l \rtimes \R^{2k}$ with the action
$\phi: \R^{2k} \rightarrow {\rm Aut}(\C^l)$ defined by
$$\phi(\bar{t}_i) ((z_1, z_2, \ldots, z_l)) = 
(e^{2 \pi \sqrt{-1}\, t_i} z_1, e^{2 \pi \sqrt{-1}\, t_i} z_2, ...,
e^{2 \pi \sqrt{-1}\, t_i} z_l),$$
where $\bar{t}_i = t_i e_i$ ($e_i$: the $i$-the unit vector in $\R^{2k})$, 
$i = 1, \ldots, 2k$. Then $\Z^{2n} \backslash G$
is a solvmanifold diffeomorphic to a torus $T^{2n}$ ($n= k+l$).
\medskip

It is easily seen that $M$ (in Example 4) is a finite quotient of a complex torus, and
has a structure of a complex torus bundle over a complex torus. The following
theorem asserts that any solvmanifolds which admit K\"ahler structures must belong
to the class of K\"ahlerian solvmanifolds defined in Example 4.
Note that a solvmanifold of completely solvable type does not belong to this class unless
it is a complex torus:
\medskip

\noindent {\bfseries Theorem 2 {\rm (\cite{H4})}.} {\em A compact solvmanifold admits a
K\"ahler structure if and only if it is a finite quotient of a complex torus which
has a structure of a complex torus bundle over a complex torus. In particular,
a compact solvmanifold of completely solvable type has a K\"ahler structure
if and only if it is a complex torus.}
\medskip

The proof of the theorem is based on a result of Aparura and Nori \cite{AN} that
a polycyclic K\"ahler group has a nilpotent subgroup of finite index, together with
our previous results: a result of the papers \cite{BG1, H1} that the only K\"ahlerian
nilmanifold is a complex torus, and a partial result on K\"ahlerian solvmanifolds
in the paper \cite{H3}. For the details of the proof and the related topics,
we refer to the paper \cite{H4}.
\medskip

We give an example of a six-dimensional solvmanifold
which admits a psuedo-K\"ahler structures but no K\"ahler structures.
\medskip

\noindent {\bfseries Example 5 ({\rm De And\'res et al. \cite{AFLM})}.}
Let $\Lambda = \Z^4 \rtimes \Z$, with the action
$\phi: \Z \rightarrow {\rm Aut}(\Z^4)$ defined by $\phi(1)=A \oplus A \in {\rm Aut}(\Z^4)$,
where $A \in {\rm SL}(2, \Z)$ has two positive eigenvalues $a_1, a_2 \;(a_1 \not=1,
a_1a_2=1)$. Let $u=(u_1, u_2), v=(v_1, v_2)$ be the eigenvectors of $a_1, a_2$
respectively. We have a basis
$\{(u_1, v_1, 0, 0), (u_2, v_2, 0, 0), (0, 0, u_1, v_1), (0, 0, u_2, v_2)\}$ of $\R^4$,
which defines an abelian lattice $\Z^4$ of $\R^4$. We define a solvable Lie group
$H=\R^4 \rtimes \R$ with the action $\bar{\phi}: \R \rightarrow {\rm Aut}(\R^4)$ defined
by $$\bar{\phi}(t)(x_1, x_2, y_1, y_2) = (e^{t{\rm log}a_1} x_1, e^{t{\rm log}a_2} x_2,
e^{t{\rm log}a_1} y_1, e^{t{\rm log}a_2} y_2),$$ 
which is a canonical extension of $\phi$. Then $\Gamma=\Lambda \times \Z$ is a lattice
of a solvable Lie group $G=H \times \R$, and $M = \Gamma \backslash G$ is a
six-dimensional solvmanifold.
\smallskip

The Lie algebra  $\mathfrak g$ of $G$ is expressed as having a basis
$\{X_1,X_2,Y_1,Y_2, Z, W\}$  with the bracket multiplication:
$$[X_1, Z] = X_1, \;[X_2, Z] = -X_2, \;[Y_1, Z] = Y_1, \;[Y_2, Z] = -Y_2,$$
and all other brackets vanishing. Let $\alpha_1, \alpha_2, \beta_1, \beta_2, \gamma, \eta$
be the corresponding Maurer-Cartan forms (left-invariant 1-forms), then we have
$$d\alpha_1 = \gamma \wedge \alpha_1, \;d\alpha_2 = -\gamma \wedge \alpha_2,\;
d\beta_1 = \gamma \wedge \beta_1,\;
d\beta_2 = -\gamma \wedge \beta_2, \;d\gamma = d \eta = 0.$$

Let $\omega = \alpha_1\wedge\alpha_2 + \beta_1\wedge\beta_2 + \gamma\wedge\eta$,
then $\omega$ defines a left-invariant symplectic form on $M$. A complex
structure $J$ on $\mathfrak g$ is defined by
$$JX_1 = Y_1, JY_1 = -X_1, JX_2 = Y_2, JY_2 = - X_2, JZ = W, JW = -Z.$$
It is easy to check that the Nijenhuis tensor vanishes for $J$, and that the pair
$(\omega, J)$ defines a psuedo-K\"ahler structure on $M$. Since $G$ is
of completely solvable type, we can apply a result of Hattori \cite{HT} that
the De Rham cohomology ring of $M$ is isomorphic to the cohomology ring of
$\mathfrak g$. It is now not hard to check that with respect to $(\omega, J)$
the following  ``K\"ahler conditions" hold for $M$ (while $M$ never admit K\"ahler
structures due to Theorem 1):
(1) the Hard Lefschetz condition, (2) the Betti number $b_{2k-1}$ are even,
(3) the De Rham complex is formal,
(4) the Fr\"ohlicher spectral sequence degenerates at $E^1$.
For more details we refer to the papers \cite{AFLM, AN}. 

\section{Four-dimensional compact homogeneous
manifolds with complex structures}

In this section we determine all the complex surfaces which are diffeomorphic
to four-dimensional compact homogeneous manifolds. 
\smallskip

It is known (due to V. V. Gorbatsevich \cite{G, GO}) that
a four-dimensional compact homogeneous manifold is diffeomorphic to one
of the following types: (1) $\prod S^{k_i}$ (up to finite quotient), where
$k_i \ge 1$ with $\sum {k_i} = 4$, (2) $\C {\rm P}^2$,
(3) Solvmanifold, (4) $S^1 \times \Gamma \backslash \widetilde{\rm SL}_2(\R)$,
where $\widetilde{\rm SL}_2(\R)$ is the universal covering of
${\rm SL}_2(\R)$ and $\Gamma$ is a lattice of $\widetilde{\rm SL}_2(\R)$.
We can determine, based on the above result, the complete list of
complex surfaces with diffeomorphism type of four-dimensional compact
homogeneous manifolds:

\renewcommand{\multirowsetup}{\centering}
$$
\begin{array}{c|c|c|c}
\hline
S & b_1 & {\rm Complex\: Structure} & \kappa \\
\hline
S^2 \times T^2 & 2 & {\rm Ruled\: Surface\: of\: genus\: 1} &
\multirow{5}{7mm}{$\left. \rule{0mm}{12mm} \right\} {-\infty}$} \\ \cline{1-3}
S^1 \times_{\Z_m} S^3/H & 1 & {\rm Hopf\: Surface} & \\ \cline{1-3}
S^2 \times S^2 & 0 & {\rm Hirzebruch\: Surface\: of\: even\: type} & \\ \cline{1-3}
\C {\rm P}^2 & 0 & {\rm Complex\: Projective\: Space} & \\ \cline{1-3}
\multirow{5}{26mm}{\hbox{Solvmanifold$\left. \rule{0mm}{12mm} \right\{$}} & 1 &
{\rm Inoue\: Surface} & \\ \cline{2-4}
& 4 & {\rm Complex\: Torus} & 
\multirow{4}{7mm}{$\left. \rule{0mm}{9mm} \right\} 0$} \\ \cline{2-3}
& 3 & {\rm Primary\: Kodaira\: Surface} & \\ \cline{2-3}
& 2 & {\rm Hyperelliptic\: Surface} & \\ \cline{2-3}
& 1 & {\rm Secondary\: Kodaira\: Surface}&  \\ \hline
S^1 \times \Gamma \backslash \widetilde{\rm SL}_2(\R)& {\rm odd} &
{\rm Properly\: Elliptic\: Surface} & 1 \\
\hline
\end{array}
$$
\noindent where $\kappa$ is the Kodaira dimension of $S$, $\rm H$ is a finite
subgroup of ${\rm SU}(2)$ acting freely on $S^3$.
\medskip

It is well known (due to A. Borel and J. P. Serre) that $S^4$ has no almost
complex structure. 
For the case of $S^2 \times T^2$, it is known
(due to T. Suwa \cite{SU}) that a complex surface is diffeomorphic to a
$S^2$-bundle over $T^2$ if and only if it is a ruled surface of genus $1$.
We can see this also from the recent result
(due to R. Friedman and Z. B. Qin \cite{FQ}) that the Kodaira dimension
of complex algebraic surface is invariant up to diffeomorphism.
To be more precise, there exist two diffeomorphism types of ruled
surfaces of genus $1$: the trivial one and the non-trivial one (which correspond
to two diffeomorphism types of $S^2$-bundles over $T^2$), and the
latter is not of homogeneous space form (see \cite{SU}).
\smallskip

For the case of $S^1 \times S^3$, K.~Kodaira showed \cite{K2}
that a complex surface diffeomorphic to a finite quotient of $S^1 \times S^3$ is a
Hopf surface. Generally a Hopf surface is diffeomorphic to a fiber bundle
over $S^1$ with fiber $S^3/{\rm U}$, defined by the action
$\rho: \pi_1(S^1) \rightarrow N_ {{\rm U}(2)}({\rm U})$
with $\rho(1)$ being cyclic of order $m$, where $\rm U$ is a finite subgroup of
${\rm U}(2)$ acting freely on $S^3$: that is, $S = S^1 \times_{\Z_m} S^3/U$
(see \cite{H2}).
We can see that a Hopf surface is of homogeneous space form if and only if
$U$ is a finite subgroup of ${\rm SU}(2)$. Let $G = {\rm SU}(2) \times S^1$,
which is a compact Lie group structure on $S^3 \times S^1$.
Take a finite subgroup $\Delta = H \rtimes \Z_m$ of $G$,
where $H$ is a finite subgroup of ${\rm SU}(2)$, $\Z_m$ is a finite
cyclic subgroup of $G$ generated by $c$:
$$c= (\tau, \xi),
\tau = \left(
\begin{array}{cc}
\xi^{-1} & 0\\
0 & \xi
\end{array}
\right), \xi^m = 1,$$
and $\tau$ belongs to $N_{{\rm SU}(2)}(H)$.
$S$ is a fiber bundle over $S^1$ with fiber $S^3/H$, which has a canonical
complex structure, defining a Hopf surface. It should be noted that if $\tau$
does not belong to $H$ and $m \ge2$, then $S$ is a non-trivial bundle.
Conversely, given a Hopf surface $S$ with fiber $S^3/H$, defined by the action
$\rho$, we can assume that $\rho(1)$ is a diagonal matrix, all of which
entries are $m$-th roots of $1$ (see \cite{H2}). Then we can see that $\rho(1)$, which belongs
to $N_{{\rm U}(2)}(H)$, actually belong to $N_{ {\rm SU}(2)}(H)$.
Hence $S$ is diffeomorphic to the one constructed above.
\smallskip

For the case of $S^2 \times S^2$, it was shown
(due to Z. B. Qin \cite{Q}) that a complex surface diffeomorphic to
$S^2 \times S^2$ must be a Hirzebruch surface of even type, which is
by definition a ruled surface of genus $0$ with diffeomorphism type
$S^2 \times S^2$. As is well known, there exist two  diffeomorphism types of
ruled surfaces of genus $0$: $S^2 \times S^2$ and $\C {\rm P}^2 \#
\overline{\C {\rm P}^2}$
(which correspond to two diffeomorphism types of $S^2$-bundles over $S^2$).
A Hirzebruch surface of odd type is the surface of the latter type. We can see that
no non-trivial finite quotient of $S^2 \times S^2$ has complex structure.
It is well known (due to S. -T. Yau) that $\C {\rm P}^2$ can have only the
standard complex structure.
We have studied in detail the case of solvmanifolds in this paper.
The complex surfaces with diffeomorphism type of solvmanifolds are Inoue
surfaces for $\kappa = -\infty$ and all of those with $c_2=0$ for $\kappa = 0$.
For the case of $S^1 \times \Gamma \backslash \widetilde{\rm SL}_2(\R)$,
C.~T.~C.~Wall showed \cite{W1} that it admits a canonical complex
structure, which define a properly elliptic surface with $b_1={\rm odd}$
and $c_2=0$; and conversely any such surface with no singular fibers
is diffeomorphic to $S^1 \times \Gamma \backslash \widetilde{\rm SL}_2(\R)$
for some lattice $\Gamma$.
\medskip

\noindent{\bfseries Acknowledgments.} The author would like to thank
V. V. Gorbatsevich for valuable comments; in particular the results in Section 8
were obtained while communicating with him by e-mail. The author appreciates
very much the referee's valuable comments which lead to many improvements
to the paper. Finally the author would like to express a sincere
gratitude to A.~Tralle for research correspondence and invitation to the conference.


\medskip

\begin{flushleft}
{\tiny
DEPARTMENT OF MATHEMATICS\\
FACULTY OF EDUCATION AND HUMAN SCIENCES\\
NIIGATA UNIVERISTY, NIIGATA, JAPAN\\
{\em E-mail address: hasegawa@ed.niigata-u.ac.jp}\\
}
\end{flushleft}

\end{document}